\documentclass[11pt]{amsart}
\usepackage{amsmath, amsthm, amsfonts, amsbsy, amssymb, upref} 
\usepackage{epsfig, graphicx}
\usepackage{pdfsync, hyperref}

\newcommand{\I}{\mathrm{i}}

\newtheorem{thm}{Theorem}

\newcommand{\ee}{\mathbb{E}}

\newcommand{\ra}{\rightarrow}
\newcommand{\rr}{\mathbb{R}}

\begin{document}
\title{Stochastic solutions of the wave equation}
\author{Sourav Chatterjee}

\address{Courant Institute of Mathematical Sciences, New York University, 251 Mercer Street, New York, NY 10012}
\email{sourav@cims.nyu.edu}
\thanks{Sourav Chatterjee's research was partially supported by the  NSF grant DMS-1005312}

\begin{abstract}
Unlike the heat equation or the Laplace equation, solutions of the wave equation on general domains have no known stochastic representation. This short note gives a simple solution to this well known problem in arbitrary dimensions. The proposed representation has several shortcomings, one of which is that it does not cover all solutions. Still, it is proof that a large class of nontrivial solutions of the wave equation in general dimensions and domains may indeed be represented stochastically.
\end{abstract}

\maketitle

\section{Introduction}
An open connected subset of $\rr^d$ is called a domain. Suppose that $D\subseteq \rr^d$ is a bounded domain with boundary $\partial D$. Let $B = (B_t)_{t\ge 0}$ be standard $d$-dimensional Brownian motion, and let $\ee_x$ denote `expectation given $B_0=x$'. Let $\tau$ denote the first time $B_t$ hits $\partial D$. If $f: \partial D \ra \rr$ is a bounded continuous function, it is well known that the function $g:D \ra\rr$ defined as
\[
g(x) := \ee_x f(B_\tau)
\]
solves the Laplace equation
\[
\Delta g(x) = 0
\]
in $D$ (where $\Delta$ denotes the $d$-dimensional Laplacian), and approaches $f$ at the boundary of $D$. This is a famous example of a stochastic representation of the solution of a partial differential equation. Another example is the heat equation
\[
\partial_t u = \frac{1}{2}\Delta u,
\]
where $u:\rr_+ \times\rr^d \ra\rr$ is a function of two variables $t\in \rr_+ = (0,\infty)$ and $x\in \rr^d$, and $\Delta$ denotes the Laplacian with respect to the $x$ variable. This problem is solved by 
\[
u(t,x) = \ee_x f(B_t),
\]
where $f:\rr^d \ra \rr$ is an arbitrary function satisfying some mild conditions. 

The goal of this note is to present a stochastic representation of certain nontrivial solutions of the wave equation
\[
\partial_t^2 u = \Delta u
\]
on a bounded domain $D$, in the spirit of the solutions of Laplace equation and the heat equation described above. Here $u:\rr_+ \times D \ra \rr$ is a function of two variables $t\in \rr_+$ and $x\in D$, and by nontriviality we mean that $u$ is not a trivial solution of the form $u(t,x)=g(x)$ for some harmonic function $g$ on $D$. 

Although the problem is well known, there are only a few papers attempting to provide stochastic representations of hyperbolic partial differential equations. Goldstein \cite{goldstein51} and Kac \cite{kac74} gave stochastic solutions of the one-dimensional telegrapher's equation. Recently, Dalang, Mueller and Tribe~\cite{dalang} gave stochastic solutions to the wave equation in dimensions 1, 2  and 3 using the solutions available from classical PDE theory, and Bakhtin and Mueller~\cite{bakhtin} used `stochastic cascades' to solve one-dimensional semilinear wave equations. The beautiful recent manuscript of Pal and Shkolnikov~\cite{pal} gives a connection between hyperbolic PDE's and intertwined diffusion processes, although it has a somewhat different approach than the aforementioned solutions of the heat and Laplace equations. 

The following theorem gives a probabilistic way of generating a large class of nontrivial solutions of wave equations on bounded domains in arbitrary dimensions. The `data' required for generating a solution is an `evolution on the boundary' in the form of a function $f(t,x)$ of two variables $t\in \rr$ and $x\in \partial D$. 
\begin{thm}\label{wavethm}
Take any $d\ge 1$ and let $D$ be a bounded open connected subset of $\rr^d$ with boundary $\partial D$. Let $f:\rr\times \partial D \ra \rr$ be any bounded measurable function. Let $B = (B_t)_{t\ge 0}$ be a $d$-dimensional standard Brownian motion, and let $\ee_x$ denote the conditional expectation given that $B_0=x$. Let $X$ be a standard Cauchy random variable and $Z$ be a standard normal random variable, independent of $B$ and independent of each other. Let $\tau$ be the minimum $t$ such  that $B_t\in \partial D$. For $t\in \rr_+$ and $x\in D$, define 
\[
u(t,x) := \ee_x f(tX + \sqrt{\tau} Z, \, B_\tau). 
\]
Then $u$ is twice continuously differentiable in all coordinates and satisfies the wave equation $\partial_t^2 u = \Delta u$ in $\rr_+\times D$. 
\end{thm}
The arbitrariness in the choice of $f$ indicates that the above procedure can be used to generate a large class of nontrivial solutions. However, the connection between $f$ and the `initial data' $u(0,x)$ and $\partial_t u(0,x)$ is somewhat complicated. It is desirable to have a representation that solves the wave equation with given values of $u(0,x)$ and $\partial_t u(0,x)$. The above solution does not fulfill that need. Still, it demonstrates that  a large class of solutions may indeed be generated stochastically. 

In general, it should not be difficult to compute $u$ by simulating Brownian motion on a computer. Sometimes explicit computations may also be possible. The following example illustrates such a case. Take $d=1$ and let $D$ be the interval $(-1,1)$. Let $f:\rr \times \{-1,1\} \ra \rr$ be the function $f(t,x) = e^x \cos t$. Note that
\begin{align*}
&\ee_x \bigl(f(tX + \sqrt{\tau} Z,\, B_\tau)\mid X, B\bigr) \\
&= \frac{e^{B_\tau}}{2}\bigl(e^{\I tX}\ee_x( e^{\I\sqrt{\tau} Z} \mid X, B) + e^{-\I tX}\ee_x( e^{-\I\sqrt{\tau} Z} \mid X, B)\bigr)\\
&= \frac{e^{B_\tau}}{2}\bigl(e^{\I tX} + e^{-\I tX}\bigr)e^{-\frac{1}{2}\tau}. 
\end{align*} 
Since $(e^{B_t - \frac{1}{2}t})_{t\ge 0}$ is a martingale, an easy application of the optional stopping theorem shows that 
\[
\ee_x(e^{B_\tau-\frac{1}{2}\tau}) = e^x. 
\]
Finally, recall that $\ee(e^{\I t X}) = e^{-|t|}$ for all $t\in \rr$. Combining the steps gives that for $t> 0$ and $x\in (-1,1)$, 
\[
u(t,x) = e^{x-t}.
\]
This, of course, is a nontrivial solution of the wave equation in $\rr_+\times (-1,1)$. It is natural to ask whether all solutions may be obtained in this manner. The simple answer is no, because the function $e^{x+t}$ is also a solution, but being unbounded,  it cannot arise according the prescription given by Theorem~\ref{wavethm}.

Can all bounded solutions be represented as in Theorem \ref{wavethm}? The answer is still no, but that is a little harder to see. Let $f:\rr\times \{-1,1\}\ra \rr$ be any bounded measurable function and let $u$ be defined as in Theorem \ref{wavethm}. Let 
\[
v(t,x) := \ee_x f(t + \sqrt{\tau}Z, \, B_\tau),
\]
so that 
\[
u(t,x) = \ee v(tX,x) = \int_{-\infty}^\infty \frac{v(tz,x)}{\pi(1+z^2)} dz = \int_{-\infty}^\infty \frac{tv(y,x)}{\pi(t^2+y^2)} dy. 
\]
Differentiating with respect to $t$ gives
\[
\partial_t u(t,x) = \int_{-\infty}^\infty \frac{(y^2-t^2)v(y,x)}{\pi(t^2+y^2)^2} dy.
\]
For any $t\ge 1$ and $y\in \rr$,
\[
\biggl|\frac{y^2-t^2}{\pi(t^2+y^2)^2}\biggr|\le \frac{1}{\pi(t^2+y^2)}\le \frac{1}{\pi(1+y^2)}. 
\]
Since $v$ is a bounded function, this shows that by the dominated convergence theorem, 
\[
\lim_{t\ra\infty}\partial_t u(t,x) = 0. 
\]
The function $u(t,x) = \cos x \cos t$ is a solution of the wave equation that does not satisfy the above property. This shows that not all bounded solutions of the wave equation are of the form given in Theorem \ref{wavethm}.  It would be interesting to characterize the full set of solutions that may be obtained using Theorem \ref{wavethm}.

\section{Proof of Theorem \ref{wavethm}}
For $t\in \rr$ and $x\in D$, define 
\[
v(t,x) := \ee_x f(t + \sqrt{\tau} Z, \, B_\tau). 
\]
Then clearly
\begin{equation}\label{uv}
u(t,x) = \ee v(tX, x).
\end{equation}
Let $W = (W_s)_{s\ge 0}$ be a one-dimensional standard Brownian motion starting at $0$, independent of $B$. Then conditional on  $B$, $W_\tau$ has the same distribution as $\sqrt{\tau} Z$. Thus,
\[
v(t,x) = \ee_x f(t+W_\tau, B_\tau). 
\]
Let $V= (V_s)_{s\ge 0}$ be a $(d+1)$-dimensional standard Brownian motion. Let $D' = \rr \times D$. Then $D'$ is a domain in $\rr^{d+1}$. Let $\sigma$ be the minimum $s$ such that $V_s\in \partial D' = \rr \times \partial D$. Take any $t\in \rr$ and $x\in D$. It is easy to see that the distribution of $V_\sigma$ given $V_0=(t,x)$ is the same as that of $(t+W_\tau, \, B_\tau)$ given $B_0=x$. Thus,
\[
v(t,x)=\ee_{(t,x)} f(V_\sigma),
\]
where $\ee_{(t,x)}$ denotes conditional expectation given $V_0=(t,x)$. By the strong Markov property of Brownian, the right hand side in the above display is a harmonic function in $D'$ (see e.g.~\cite[Theorem 3.8]{morters}). Consequently, $v$ is harmonic in $\rr \times D$, and therefore
\begin{equation}\label{veq}
\partial_t^2 v = -\Delta v, 
\end{equation}
where $\Delta$ denotes the Laplacian in the $x$ coordinates. 

For each $a> 0$, define 
\[
v_a(t,x) := e^{-at^2} v(t,x)
\]
and 
\[
u_a(t,x) := \ee v_a(tX, x). 
\]
Since $f$ is a bounded function, so is $v$. Take any $x\in D$. Then there is an $r> 0$ such that for any $t\in \rr$, the Euclidean ball of radius $r$ around $(t,x)$ is contained in $D'$. Using this and the boundedness of $v$, it follows by standard estimates for harmonic functions (see e.g.~\cite[Corollary 1.4]{bass}) that any derivative of $v$, of any order, is uniformly bounded as $t$ varies arbitrarily in $\rr$ and $x$ varies inside a small ball away from the boundary of $D$. (Recall that a harmonic function is, in particular, $C^\infty$.)

The uniform boundedness of derivatives of $v$ implies the same for the derivatives of $v_a$ for any $a> 0$. Moreover, the bounds will not depend on $a$. These boundedness properties are important ingredients for the rest of the proof; they will be collectively referred to as the `uniform boundedness of derivatives'.

Since $v$ is bounded, it follows by the dominated convergence theorem that
\begin{equation}\label{ualim}
u(t,x) = \lim_{a\ra0} u_a(t,x).
\end{equation}
A simple consequence of the uniform boundedness of derivatives is that  
\begin{equation}\label{uava}
\Delta u(t,x) = \ee \Delta v(tX, x), 
\end{equation}
and also that $\Delta u (t,x)$ is a continuous function of $t$. 
Define
\[
v_a^{(1)} := \partial_t v_a, \ \  v_a^{(2)} := \partial_t^2 v_a.
\]
Then  observe that for any $a> 0$ and $t> 0$, the uniform boundedness of derivatives and the rapidly decaying property of $e^{-at^2}$ imply that while differentiating $u_a$ with respect to $t$, the derivative may be moved inside the expectation to give
\begin{align}
\partial_t^2 u_a(t,x) &= \ee (X^2v_a^{(2)}(tX, x))\nonumber \\
&= \frac{1}{\pi}\int_{-\infty}^\infty v_a^{(2)}(tz, x) dz - \int_{-\infty}^\infty \frac{v_a^{(2)}(tz,x)}{\pi(1+z^2)} dz. \label{mainequation}
\end{align}
Again, the rapidly decaying property of $e^{-at^2}$ and the boundedness of the derivatives of $v$ imply that 
\[
\int_{-\infty}^\infty |v_a^{(2)}(tz, x)| dz <\infty,
\]
and therefore
\begin{align}
\int_{-\infty}^\infty v_a^{(2)}(tz, x) dz &= \lim_{R\ra\infty}\int_{-R}^R v_a^{(2)}(tz, x) dz\nonumber \\
&= \frac{1}{t} \lim_{R\ra\infty} (v_a^{(1)}(tR, x)- v_a^{(1)}(-tR, x))\nonumber \\
&= 0. \label{va2}
\end{align}
Equations \eqref{mainequation} and \eqref{va2} show that for any $a> 0$ and $t> 0$, 
\begin{equation}\label{e1}
\partial_t^2u_a(t,x) = - \ee v_a^{(2)}(tX, x). 
\end{equation}
A simple computation using equation \eqref{veq} shows that for all $t\in \rr$ and $x\in D$,
\begin{equation}\label{e2}
\lim_{a\ra 0} v_a^{(2)} (t,x) = - \Delta v(t,x). 
\end{equation}
By the uniform boundedness of derivatives, the dominated convergence theorem, and equations \eqref{uava}, \eqref{e1} and \eqref{e2}, \begin{equation}\label{ulim}
\lim_{a\ra 0} \partial_t^2 u_a(t,x) =  \ee\Delta v(tX,x) = \Delta u(t,x). 
\end{equation}
Next, take any $0<t_0 < t$. Then
\begin{align*}
u_a(t,x) - u_a(t_0,x) - (t-t_0)\partial_t u_a(t_0,x)  &= \int_{t_0}^t \int_{t_0}^s \partial_t^2 u_a(w,x) \,dw \,ds. 
\end{align*}
By \eqref{ualim}, 
\[
\lim_{a\ra 0} u_a(t,x)=u(t,x), \ \ \lim_{a\ra0} u_a(t_0,x) = u(t_0,x).
\]
By \eqref{ulim} and the boundedness of derivatives, 
\[
\lim_{a\ra0}  \int_{t_0}^t \int_{t_0}^s \partial_t^2 u_a(w,x) \,dw \,ds =  \int_{t_0}^t \int_{t_0}^s \Delta u(w,x) \,dw \,ds. 
\]
Thus, $\lim_{a\ra0}\partial_t u_a(t_0,x)$ exists. Let this limit be denoted by $C(t_0,x)$. Then by the last three displays, 
\[
u(t,x) - u(t_0,x) - (t-t_0)C(t_0,x)  = \int_{t_0}^t \int_{t_0}^s \Delta u(w,x) \,dw \,ds.
\]
The above equation and the continuity of $\Delta u(t,x)$ in $t$  shows that $u(t,x)$ is twice continuously differentiable in $t$ and
\[
\partial_t^2 u = \Delta u,
\]
as required. This completes the proof. 
\vskip.3in
\noindent {\bf Acknowledgments.} The authors thanks Misha Shkolnikov for pointing out the reference \cite{goldstein51}.

\end{document}